\documentclass[12pt]{amsart}
\usepackage{latexsym, amsmath, amssymb}
\usepackage{epsfig}
\usepackage{amscd}
\usepackage{pstricks,pst-node,cite}
\newtheorem{thm}{Theorem}[section]
\newtheorem{lem}[thm]{Lemma}

\newtheorem{exmp}[thm]{Example}

\newtheorem{conj}[thm]{Conjecture}
\theoremstyle{remark}

\newcommand{\Inv}{\mbox{\rm Inv}}

\newcounter{fignum}
\ifx\blackandwhite\undefined
 \psset{linecolor=blue}\newrgbcolor{darkred}{.9 0 0}\newrgbcolor{emgreen}{0 .9 0}
 \newrgbcolor{pup}{.7  0 .9}\newrgbcolor{xxx}{.8 .8 .8} \else
   \psset{linecolor=black}\newgray{darkred}{.7}\newgray{emgreen}{0.3}\newgray{pup}{.5}\newgray{xxx}{.5}
\fi

\psset{linewidth=.5pt,dash=3pt 3pt,doublesep=.05,dotsize=1pt 5}
\SpecialCoor

\newcommand{\mda}{\lput{:U}{\pspicture(0,0)(0,0)
\psline[linecolor=darkred,arrows=->,arrowscale=1.7](3.2pt,0)(3.4pt,0)\endpspicture}}
\newcommand{\mdb}{\lput{:U}{\pspicture(0,0)(0,0)
\psline[linecolor=emgreen,arrows=->,arrowscale=1.7](3.2pt,0)(3.4pt,0)\endpspicture}}
\newcommand{\middlearrow}{\lput{:U}{\pspicture(0,0)(0,0)
\psline[linecolor=blue,arrows=->,arrowscale=1.7](3.2pt,0)(3.4pt,0)\endpspicture}}

\begin{document}

\author{Dongseok KIM}
\address{Department of Mathematics \\ Kyungpook National University \\ Taegu 702-201 Korea}
\email{dongseok@mail.knu.ac.kr}
\thanks{}
\author{Jaeun Lee}
\address{Department of Mathematics\\ Yeungnam University\\ Kyongsan, 712-749, Korea }
\email{julee@yu.ac.kr}
\thanks{The first author was supported in part by KRF Grant M02-2004-000-20044-0,
the second author was supported in part by Com$^2$Mac-KOSEF}
\title[The quantum $\mathfrak{sl}(3)$ invariants of cubic bipartite planar graphs]
{The quantum $\mathfrak{sl}(3)$ invariants of cubic bipartite planar
graphs}
\begin{abstract}
Temperley-Lieb algebras have been generalized to
$\mathfrak{sl}(3,\mathbb{C})$ web spaces. Since a cubic bipartite
planar graph with suitable directions on edges is a web, the quantum
$\mathfrak{sl}(3)$ invariants naturally extend to all cubic
bipartite planar graphs. First we completely classify them as a
connected sum of primes webs. We also provide a method to find all
prime webs and exhibit all prime webs up to $20$ vertices. Using
quantum $\mathfrak{sl}(3)$ invariants, we provide a criterion which
determine the symmetry of graphs.
\end{abstract}

\maketitle
\section{introduction}

Triangulations are maximal simple planar graphs where we can not add
an edge without destroying the planarity. Every triangulations are
$3$-connected, thus it has a unique embedding into the two
dimensional sphere \cite{whitney}. The dual of a triangulation is a
$3$-connected cubic planar graphs. If a cubic planar graph is not
$3$-connected, it could have several nonequivalent embeddings.
Decomposing these cubic planar graphs by the connectivity has
powerful applications such as to count all cubic planar graphs
\cite{Walshcounting}. The Barnette's conjecture, every $3$-connected
cubic bipartite planar graphs are Hamiltonian \cite{barnette},
brought a lot of attentions to cubic bipartite planar graphs
\cite{goodey}.

\begin{figure}
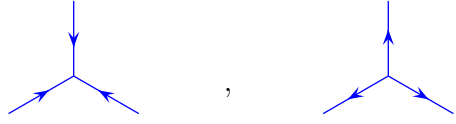

$$
\pspicture[.4](-1,-1)(1,1) \pcline(0,1)(0,0)\middlearrow
\pcline(-.866,-.5)(0,0)\middlearrow
\pcline(.866,-.5)(0,0)\middlearrow
\endpspicture
\hskip 1cm , \hskip 1cm \pspicture[.4](-1,-1)(1,1)
\pcline(0,0)(0,1)\middlearrow \pcline(0,0)(-.866,-.5)\middlearrow
\pcline(0,0)(.866,-.5)\middlearrow
\endpspicture
$$
\caption{Generators of the quantum $\mathfrak{sl}(3)$ web space.}
\label{generator}
\end{figure}

An unexpected relation between the representation theory and
bipartite cubic planar graphs has disclosed as follows. After the
discovery of the Jones polynomial~\cite{Jonespoly}, its
generalizations have been studied in many different ways. One of
successful generalizations is to use the representation theory of
complex simple Lie algebras from the original work of Reshetikhin
and Turaev \cite{RT1}. In particular, Kuperberg generalized
Temperley-Lieb algebras, which corresponds to the invariant subspace
of a tensor product of the vector representation of
$\mathfrak{sl}(2)$, to web spaces of simple Lie algebras of rank
$2$, $\mathfrak{sl}(3)$, $\mathfrak{sp}(4)$ and $G_2$
\cite{Kuperbergspiders}. All webs in the representation theory of
$\mathfrak{sl}(3)$ are generated by the webs in Figure
\ref{generator} with a complete set of the relations presented in
Figure~\ref{relations} where the quantum integers are defined as
$$
[n]=\frac{q^{\frac{n}{2}}-q^{-\frac{n}{2}}}{q^{\frac{1}{2}}-q^{-\frac{1}{2}}}.
$$
Roughly speaking, a \emph{basis web} is a web for which we can not
apply any relation in Figure~\ref{relations}. For precise
definition, we refer to~\cite{Kuperbergspiders}. There are many
interesting results on the quantum $\mathfrak{sl}(3)$ invariants
\cite{chbili, Khovanovsl3, Khovanovcolored, Dongseokthesis,
Leintegral, Lickorishsu, MOYHomfly, OYquantum, Sokolov}.

The $\mathfrak{sl}(3)$ webs are directed cubic bipartite planar
graphs together with circles (without vertices, thus different from
loops) where the direction of the edges is from one set to the other
set in the bipartition. From a given cubic bipartite planar graph,
one can find its quantum $\mathfrak{sl}(3)$ invariant as follows.
Once we direct the edges by one of two possible orientations as we
described, it can be considered a web of empty boundary. Since the
dimension of the web space of webs of empty boundary is one over
$\mathbb{C}[q^{\frac12},q^{-\frac12}]$, we obtain a Laurent
polynomial. In Lemma~\ref{direction}, we prove the quantum
$\mathfrak{sl}(3)$ invariant does not depend on the choice of
directions. Therefore, the quantum $\mathfrak{sl}(3)$ invariant
naturally extends to any cubic bipartite planar graph $G$, let us
denote it by $P_G(q)$. If a cubic bipartite planar graph $G$ is a
disjoint union of two cubic bipartite planar graphs $G_1$ and $G_2$,
then $P_G = P_{G_1}P_{G_2}$. Thus, we assume all abstract cubic
bipartite planar graphs are connected. Then we find the following
classification theorem.

\begin{thm}
Let $G$ be a connected cubic bipartite planar graph. Then there
exist $3$-connected cubic bipartite planar graphs $G_1, G_2, \ldots,
G_k$ such that
$$G = G_1 \# G_2 \# \ldots \# G_k,$$
where the $\#$ operation is defined in Figure \ref{connectedsum}.
Moreover, the decomposition does not depend the choice of the planar
imbedding of $G$ up to a reflection of $G_i$'s on the two
dimensional sphere, thus the decomposition is unique up to the
reflections of $G_i$. For the quantum $\mathfrak{sl}(3)$ graph
invariant of $G$, we find
$$[3]^{k-1}P_G(q)= (-[2])^{l}P_{G_1}(q) P_{G_2}(q) \ldots P_{G_k}(q),$$
where $l$ is the number of times we use relation~\ref{a2defr22}
shown in Figure~\ref{relations} in the process of the
decomposition. \label{mainthm11}
\end{thm}

\begin{figure}
\begin{eqnarray}
\pspicture[.45](-.6,-.5)(.6,.5)
\pscircle(0,0){.4}\psline[arrows=->,arrowscale=1.5](.1,.4)(.11,.4)
\endpspicture
& = & [3] \label{a2defr21} \\
 \pspicture[.45](-1.5,-.8)(1.5,.8)
\qline(-1.2,0)(-.4,0)\psline[arrowscale=1.5]{->}(-.7,0)(-.9,0)
\pccurve[angleA=90,angleB=90,nodesep=1pt](-.4,0)(.4,0)\middlearrow
\pccurve[angleA=-90,angleB=-90,nodesep=1pt](-.4,0)(.4,0)\middlearrow
\qline(1.2,0)(.4,0)\psline[arrowscale=1.5]{->}(.9,0)(.7,0)
\endpspicture
& = & - [2] \pspicture[.45](-.8,-.8)(.8,.8) \qline(.6,0)(-.6,0)
\psline[arrowscale=1.5]{->}(.1,0)(-.1,0)
\endpspicture  \label{a2defr22} \\
\pspicture[.45](-1.1,-1.1)(1.1,1.1) \qline(1,1)(.5,.5)
\qline(-1,1)(-.5,.5) \qline(-1,-1)(-.5,-.5) \qline(1,-1)(.5,-.5)
\qline(.5,.5)(-.5,.5) \qline(-.5,.5)(-.5,-.5)
\qline(-.5,-.5)(.5,-.5) \qline(.5,-.5)(.5,.5)
\psline[arrowscale=1.5]{<-}(.1,.5)(-.1,.5)
\psline[arrowscale=1.5]{<-}(-.1,-.5)(.1,-.5)
\psline[arrowscale=1.5]{<-}(.5,.1)(.5,-.1)
\psline[arrowscale=1.5]{<-}(-.5,-.1)(-.5,.1)
\psline[arrowscale=1.5]{->}(.85,.85)(.65,.65)
\psline[arrowscale=1.5]{<-}(-.85,.85)(-.65,.65)
\psline[arrowscale=1.5]{->}(-.85,-.85)(-.65,-.65)
\psline[arrowscale=1.5]{<-}(.85,-.85)(.65,-.65)
\endpspicture
&= &  \pspicture[.45](-1.1,-1.1)(1.1,1.1)
\pnode(1;45){a1}\pnode(1;135){a2}\pnode(1;225){a3}\pnode(1;315){a4}
\nccurve[angleA=225,angleB=315]{a1}{a2}\middlearrow
\nccurve[angleA=45,angleB=135]{a3}{a4}\middlearrow
\endpspicture + \pspicture[.45](-1.1,-1.1)(1.1,1.1)
\pnode(1;45){a1}\pnode(1;135){a2}\pnode(1;225){a3}\pnode(1;315){a4}
\nccurve[angleA=225,angleB=135]{a1}{a4}\middlearrow
\nccurve[angleA=45,angleB=315]{a3}{a2}\middlearrow
\endpspicture
\label{a2defr23}
\end{eqnarray}
\caption{Complete relations of the quantum $\mathfrak{sl}(3)$ web
space.} \label{relations}
\end{figure}

The outline of this paper is as follows. In section~\ref{graph} we
first classify cubic bipartite planar graphs and show $P_G$ can be
computed by the decomposition. Using the representation theory, we
provide a method to find all prime webs and exhibit all prime webs
up to $20$ vertices in section~\ref{primesec}. We provide a
criterion which determine the symmetry of cubic bipartite planar
graphs and we discuss a few problems in section~\ref{discussion}.

\section{A classification of cubic bipartite planar graphs}
\label{graph}

First we prove the quantum $\mathfrak{sl}(3)$ invariants can be
defined for (undirected) cubic bipartite planar graphs. For terms
and notations for graph theory, we refer to \cite{GT1}.

\begin{lem}\label{direction}
$P_G(q)$ does not depend on the direction of the edges.
\end{lem}
\begin{proof}
From a given connected cubic bipartite planar graph, there are
exactly two possible ways to direct the edges to get a web and one
can be obtained from the other by reversing the directions of all
edges. The first two relations \ref{a2defr21}, \ref{a2defr22} shown
in Figure~\ref{relations} can be applied exactly same way if we
reverse the orientations of all edges. For the relation
\ref{a2defr23}, it does not change how the rectangle splits but it
only changes the direction of the edges. Therefore, $P_G$ does not
depend on the direction of the edges.
\end{proof}

By comparing these three relations for the mirror image, we find the
following lemma.

\begin{lem}\label{mirror}
Let $G$ be a cubic bipartite planar graph and let $\tilde{G}$ be the
mirror image of $G$. Then,
$$P_{G}(q)=P_{\tilde{G}}(q).$$
\end{lem}

\begin{figure}
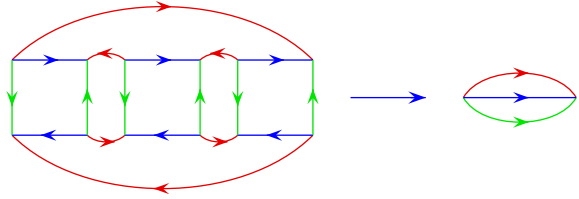

$$
\pspicture[.45](-2.2,-.5)(5.7,2.1)
\psline(2.5,1)(3.5,1)\psline[arrowscale=1.9]{->}(3.3,1)(3.5,1)
\pccurve[angleA=45,angleB=135,linecolor=darkred](-2,1.5)(2,1.5)\mda
\pccurve[angleA=135,angleB=45,linecolor=darkred](1,1.5)(.5,1.5)\mda
\pccurve[angleA=135,angleB=45,linecolor=darkred](-.5,1.5)(-1,1.5)\mda
\pccurve[angleA=-45,angleB=-135,linecolor=darkred](.5,.5)(1,.5)\mda
\pccurve[angleA=-45,angleB=-135,linecolor=darkred](-1,.5)(-.5,.5)\mda
\pccurve[angleA=-135,angleB=-45,linecolor=darkred](2,.5)(-2,.5)\mda
\pcline(1,1.5)(2,1.5)\middlearrow \pcline(2,.5)(1,.5)\middlearrow
\pcline(-.5,1.5)(.5,1.5)\middlearrow
\pcline(.5,.5)(-.5,.5)\middlearrow
\pcline(-2,1.5)(-1,1.5)\middlearrow
\pcline(-1,.5)(-2,.5)\middlearrow
\pcline[linecolor=emgreen](-.5,1.5)(-.5,.5)\mdb
\pcline[linecolor=emgreen](.5,.5)(.5,1.5)\mdb
\pcline[linecolor=emgreen](1,1.5)(1,.5)\mdb
\pcline[linecolor=emgreen](2,.5)(2,1.5)\mdb
\pcline[linecolor=emgreen](-2,1.5)(-2,.5)\mdb
\pcline[linecolor=emgreen](-1,.5)(-1,1.5)\mdb
\rput(4,1){\rnode{c1}{$$}} \rput(5.5,1){\rnode{c2}{$$}}
\nccurve[angleA=60,angleB=120,linecolor=darkred]{c1}{c2}\mda
\nccurve[angleA=-60,angleB=-120,linecolor=emgreen]{c1}{c2}\mdb
\ncline{c1}{c2}\middlearrow
\endpspicture
$$
\caption{A graph covering structure of a prime web.}
\label{covering}
\end{figure}

Let $G$ be a web of empty boundary. As we mentioned, $G$ is a
cubic bipartite directed planar graph where the directions of
edges are from one set to the other set in the definition of the
bipartition. By Lemma \ref{direction} we assume one of these
directions on the edges when we say a cubic bipartite planar
graph. Because of the directions of edges, there does not exist a
loop (it is different from the circle without a vertex). Since we
can remove all multiple edges by relation \ref{a2defr22} presented
in Figure~\ref{relations}, we assume $G$ is simple. The
\emph{connectivity} $\kappa(G)$ of a graph $G$ is the minimum
number of vertices whose removal results in a disconnected graph
or a single vertex. A graph is \emph{n-connected} if the
connectivity of $G$ is $n$ or greater. A graph $G$ is
\emph{k-edge-connected} if the removal of fewer than $k$ edges
from $G$ still leaves a connected graph, we denote the edge
connectivity of $G$ by $\kappa'(G)$. For an $r$-regular graph $G$,
$\kappa(G)\le \kappa'(G)\le r$, hence if $\kappa(G)=r$, then
$\kappa(G)=\kappa'(G)=r$. For connected trivalent graphs, these
two connectivities are the same with one exception $\Theta$, the
base graph of the covering in Figure~\ref{covering}. Since it is
not simple and we are going to deal with only simple graphs, now
we can assume $\kappa(G) = \kappa'(G)$. Since it is cubic, its
connectivity is less than or equal to $3$. If it is $3$-connected,
by a celebrated theorem by Whitney, there is only one way to imbed
$G$ into the plane, up to a reflection \cite{whitney}. A cubic
bipartite planar graphs is \emph{prime} if it is $3$-connected.
Conventionally, we assume the circle without a vertex is not
prime. If $G$ is not $3$-connected, then we want to decompose $G$
into a smaller pieces of prime graphs.

\begin{figure}
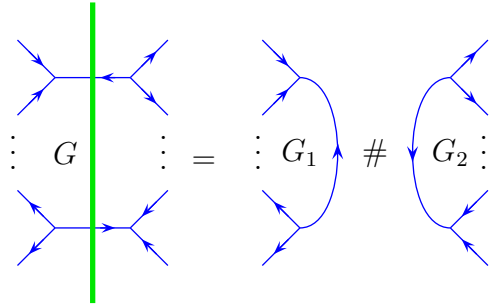

$$
\pspicture[.45](-1.2,-2.1)(1.2,2.1) \rput[r](-.2,0){$G$}
\rput[r](1,.2){$.$} \rput[r](1,0){$.$} \rput[r](1,-.2){$.$}
\rput[r](-1,.2){$.$} \rput[r](-1,0){$.$} \rput[r](-1,-.2){$.$}
\psline(-.5,1)(.5,1)\psline[arrowscale=1.5]{->}(.3,1)(.1,1)
\psline(-.5,1)(-1,1.5)\psline[arrowscale=1.5]{->}(-.85,1.35)(-.65,1.15)
\psline(-.5,1)(-1,.5)\psline[arrowscale=1.5]{->}(-.85,.65)(-.65,.85)
\psline(.5,1)(1,1.5)\psline[arrowscale=1.5]{->}(.65,1.15)(.85,1.35)
\psline(.5,1)(1,.5)\psline[arrowscale=1.5]{->}(.65,.85)(.85,.65)
\psline(-.5,-1)(.5,-1)\psline[arrowscale=1.5]{<-}(.3,-1)(.1,-1)
\psline(-.5,-1)(-1,-1.5)\psline[arrowscale=1.5]{<-}(-.85,-1.35)(-.65,-1.15)
\psline(-.5,-1)(-1,-.5)\psline[arrowscale=1.5]{<-}(-.85,-.65)(-.65,-.85)
\psline(.5,-1)(1,-1.5)\psline[arrowscale=1.5]{<-}(.65,-1.15)(.85,-1.35)
\psline(.5,-1)(1,-.5)\psline[arrowscale=1.5]{<-}(.65,-.85)(.85,-.65)
\psline[linecolor=emgreen, linewidth=2pt](0,2)(0,-2)
\endpspicture
= \pspicture[.45](-2,-2.1)(2,2.1) \rput(0,0){$\#$}
\rput(-1,0){$G_1$} \rput(1,0){$G_2$} \rput[r](1.5,.2){$.$}
\rput[r](1.5,0){$.$} \rput[r](1.5,-.2){$.$} \rput[r](-1.5,.2){$.$}
\rput[r](-1.5,0){$.$} \rput[r](-1.5,-.2){$.$}
\rput(1,1){\rnode{a1}{$$}} \rput(-1,1){\rnode{a2}{$$}}
\rput(-1,-1){\rnode{a3}{$$}} \rput(1,-1){\rnode{a4}{$$}}
\nccurve[angleA=180,angleB=180]{a1}{a4}\middlearrow
\nccurve[angleA=0,angleB=0]{a3}{a2}\middlearrow
\psline(-1,1)(-1.5,1.5)\psline[arrowscale=1.5]{->}(-1.35,1.35)(-1.15,1.15)
\psline(-1,1)(-1.5,.5)\psline[arrowscale=1.5]{->}(-1.35,.65)(-1.15,.85)
\psline(1,1)(1.5,1.5)\psline[arrowscale=1.5]{->}(1.15,1.15)(1.35,1.35)
\psline(1,1)(1.5,.5)\psline[arrowscale=1.5]{->}(1.15,.85)(1.35,.65)
\psline(-1,-1)(-1.5,-1.5)\psline[arrowscale=1.5]{<-}(-1.35,-1.35)(-1.15,-1.15)
\psline(-1,-1)(-1.5,-.5)\psline[arrowscale=1.5]{<-}(-1.35,-.65)(-1.15,-.85)
\psline(1,-1)(1.5,-1.5)\psline[arrowscale=1.5]{<-}(1.15,-1.15)(1.35,-1.35)
\psline(1,-1)(1.5,-.5)\psline[arrowscale=1.5]{<-}(1.15,-.85)(1.35,-.65)
\endpspicture
$$
\caption{A connected sum decomposition of a $2$-connected cubic
bipartite planar graph $G$ into $G_1\# G_2$.} \label{connectedsum}
\end{figure}

\begin{lem}
Let $G$ be a connected cubic bipartite planar graph. Then $G$ is
either $3$-connected or $2$-connected. Moreover, if it is
$2$-connected, then two edges, whose removal results two disconnect
graphs, intersect the separating circle, presented by a thick gray
line, in alternating directions as in Figure~\ref{connectedsum}.
\end{lem}

\begin{proof}
One can prove it by using a method in graph theory. But if we use
the representation theory, one can easily prove the result. In the
language of representation theory, the lemma can be restated that
there does not exist a cut circle of the weights $\lambda_1,
\lambda_2, 2\lambda_1$ or $2\lambda_2$. Since
 $\mathrm{dim}(\mathrm{Inv}( V_{\lambda_i}))\neq 0$,
$G$ is not $1$-connected. Since $(V_{\lambda_1})^*\cong
V_{\lambda_2}$, $(V_{\lambda_2})^*\cong V_{\lambda_1}$, by a simple
application of Schur's lemma, we find

$$\mathrm{dim}(\mathrm{Inv}( V_{\lambda_i}\otimes (V_{\lambda_j})^{*})) = \left\{
\begin{array}{cl} 1 & ~~\mathrm{if}~~ i=j, \\
0 & ~~\mathrm{if}~~ i\neq j. \end{array}\right.
$$
Moreover, if $G$ is $2$-connected, then the weight of the cut circle
must be $\lambda_1+\lambda_2$.
\end{proof}

\subsection{Proof of Theorem \ref{mainthm11}}

If $G$ is $2$-connected, then its general shape is given in the
left hand side of the equality in Figure~\ref{connectedsum}. Then
we naturally define a connected sum of two graphs (depend on the
choice of the edges). The resulting graphs $G_1$ and $G_2$ are
obviously cubic bipartite planar graphs but they may not be
simple. We apply relation \ref{a2defr22} shown in
Figure~\ref{relations} and resulting one still may not be simple
but we repeat the relation until it becomes simple. From these two
components $G_1$ and $G_2$, we repeat the process. $3$-connected
parts are not changed by relation \ref{a2defr22} shown in
Figure~\ref{relations}, thus there is no ambiguity for the order
of connected sums and the number of times we use the relation
\ref{a2defr22} does not depend neither. The finiteness of $G$
implies that the process stops and we obtain a unique prime
decomposition of $G$. The following equality immediately follows
from the decomposition process,

$$[3]^{k-1}P_G(q)= (-[2])^{l}P_{G_1}(q) P_{G_2}(q) \ldots
P_{G_k}(q).$$

One can see that this connected sum does depend on the choice of the
edges we connect two graphs. Thus, $G_1\# G_2$ is not unique in
general. Using this idea, one can construct non-isomorphic non-prime
graphs of the same quantum $\mathfrak{sl}(3)$ invariant.

\section{Prime cubic bipartite planar graphs}
\label{primesec}

\begin{figure}
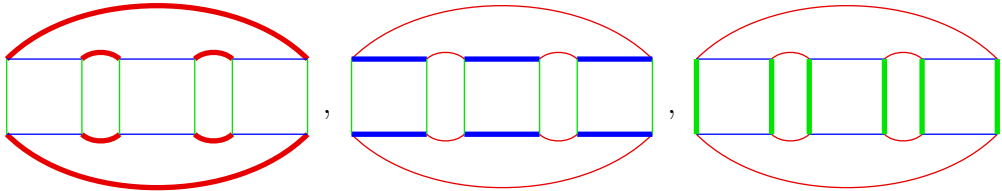

$$
\pspicture[.45](-2.2,-.3)(2.2,2.1) \rput(2,1.5){\rnode{a1}{$$}}
\pccurve[angleA=45,angleB=135,linecolor=darkred,linewidth=2pt](-2,1.5)(2,1.5)
\pccurve[angleA=135,angleB=45,linecolor=darkred,linewidth=2pt](1,1.5)(.5,1.5)
\pccurve[angleA=135,angleB=45,linecolor=darkred,linewidth=2pt](-.5,1.5)(-1,1.5)
\pccurve[angleA=-45,angleB=-135,linecolor=darkred,linewidth=2pt](.5,.5)(1,.5)
\pccurve[angleA=-45,angleB=-135,linecolor=darkred,linewidth=2pt](-1,.5)(-.5,.5)
\pccurve[angleA=-135,angleB=-45,linecolor=darkred,linewidth=2pt](2,.5)(-2,.5)
\qline(1,1.5)(2,1.5)\qline(2,.5)(1,.5) \qline(-.5,1.5)(.5,1.5)
\qline(.5,.5)(-.5,.5) \qline(-2,1.5)(-1,1.5) \qline(-1,.5)(-2,.5)
\psline[linecolor=emgreen](-.5,1.5)(-.5,.5)
\psline[linecolor=emgreen](.5,.5)(.5,1.5)
\psline[linecolor=emgreen](1,1.5)(1,.5)
\psline[linecolor=emgreen](2,.5)(2,1.5)
\psline[linecolor=emgreen](-2,1.5)(-2,.5)
\psline[linecolor=emgreen](-1,.5)(-1,1.5)
\endpspicture ,\pspicture[.45](-2.2,-.3)(2.2,2.1)
\pccurve[angleA=45,angleB=135,linecolor=darkred](-2,1.5)(2,1.5)
\pccurve[angleA=135,angleB=45,linecolor=darkred](1,1.5)(.5,1.5)
\pccurve[angleA=135,angleB=45,linecolor=darkred](-.5,1.5)(-1,1.5)
\pccurve[angleA=-45,angleB=-135,linecolor=darkred](.5,.5)(1,.5)
\pccurve[angleA=-45,angleB=-135,linecolor=darkred](-1,.5)(-.5,.5)
\pccurve[angleA=-135,angleB=-45,linecolor=darkred](2,.5)(-2,.5)
\psline[linewidth=2pt](1,1.5)(2,1.5)
\psline[linewidth=2pt](2,.5)(1,.5)
\psline[linewidth=2pt](-.5,1.5)(.5,1.5)
\psline[linewidth=2pt](.5,.5)(-.5,.5)
\psline[linewidth=2pt](-2,1.5)(-1,1.5)
\psline[linewidth=2pt](-1,.5)(-2,.5)
\psline[linecolor=emgreen](-.5,1.5)(-.5,.5)
\psline[linecolor=emgreen](.5,.5)(.5,1.5)
\psline[linecolor=emgreen](1,1.5)(1,.5)
\psline[linecolor=emgreen](2,.5)(2,1.5)
\psline[linecolor=emgreen](-2,1.5)(-2,.5)
\psline[linecolor=emgreen](-1,.5)(-1,1.5)
\endpspicture,
\pspicture[.45](-2.2,-.3)(2.2,2.1)
\pccurve[angleA=45,angleB=135,linecolor=darkred](-2,1.5)(2,1.5)
\pccurve[angleA=135,angleB=45,linecolor=darkred](1,1.5)(.5,1.5)
\pccurve[angleA=135,angleB=45,linecolor=darkred](-.5,1.5)(-1,1.5)
\pccurve[angleA=-45,angleB=-135,linecolor=darkred](.5,.5)(1,.5)
\pccurve[angleA=-45,angleB=-135,linecolor=darkred](-1,.5)(-.5,.5)
\pccurve[angleA=-135,angleB=-45,linecolor=darkred](2,.5)(-2,.5)
\qline(1,1.5)(2,1.5) \qline(2,.5)(1,.5) \qline(-.5,1.5)(.5,1.5)
\qline(.5,.5)(-.5,.5) \qline(-2,1.5)(-1,1.5) \qline(-1,.5)(-2,.5)
\psline[linecolor=emgreen,linewidth=2pt](-.5,1.5)(-.5,.5)
\psline[linecolor=emgreen,linewidth=2pt](.5,.5)(.5,1.5)
\psline[linecolor=emgreen,linewidth=2pt](1,1.5)(1,.5)
\psline[linecolor=emgreen,linewidth=2pt](2,.5)(2,1.5)
\psline[linecolor=emgreen,linewidth=2pt](-2,1.5)(-2,.5)
\psline[linecolor=emgreen,linewidth=2pt](-1,.5)(-1,1.5)
\endpspicture
$$
\caption{Three polygonal decompositions of the prime web in
Figure~\ref{covering}.} \label{decompo}
\end{figure}

For a fixed boundary, there is a systematic way to generate all webs
\cite{KKnotdual}. But we are only interested in prime webs, thus we
develop a new way to produce all prime webs in the section. As a
cubic bipartite graph, all prime webs are graph coverings of the
graph $\Theta$~\cite{KL}. Since the dual graph of a web is a planar
triangulations that the valences of all vertices are even, this dual
graph is vertex $3$-colorable, $i. e.$, cubic bipartite planar
graphs are edge $3$-colorable. In fact, the number of such colorings
of an (undirected) cubic bipartite planar graph is known
\cite{jaeger}. Thus, each polygons of a web can be alternatively
colored by just two colors naturally come from the $3$-coloring of
the dual graph. From such a coloring, we can consider a web as a
union of polygons colored by two colors and these polygons are
connected by the edges of a color which has not been used yet, we
called it a \emph{polygonal decomposition} of the web. Since there
are three colors, we find exactly three polygonal decompositions of
a web as in Figure~\ref{decompo}. In fact, the first two are
identical up to a permutation of colors.

In particular, some polygonal decompositions can be obtained from
the representation theory of the quantum $\mathfrak{sl}(2)$. First,
we fix a face which is not a polygon in a polygonal decomposition of
a web, we call it an \emph{exterior} face. For each polygon in a
polygonal decomposition, we define the \emph{level} of a polygon by
the minimum of the faces, it has to cross to reach the exterior
face. In fact, the level of a polygon is the minimal \emph{length}
of paths between vertices in the dual graph, where vertices are
corresponding to the polygon and the exterior face. A polygonal
decomposition of a web is \emph{circular} if there exists a suitable
exterior face such that all polygons in the decomposition has level
$1$. A web is \emph{circular} if at least one of polygonal
decompositions of the web is circular. All three in
Figure~\ref{decompo} are circular, so is the web.

\begin{figure}
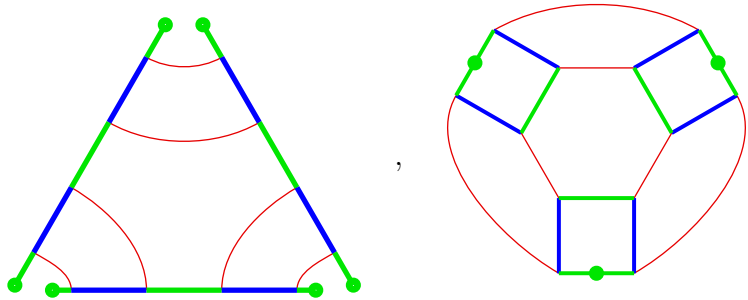

$$
\pspicture[.45](-2.8,0)(2.8,3.7)
\pccurve[angleA=-150,angleB=-30,linecolor=darkred](.5,3.098)(-.5,3.098)
\pccurve[angleA=-30,angleB=-150,linecolor=darkred](-1,2.232)(1,2.232)
\pccurve[angleA=90,angleB=-30,linecolor=darkred](-.5,0)(-1.5,1.366)
\pccurve[angleA=-30,angleB=90,linecolor=darkred](-2,.5)(-1.5,0)
\pccurve[angleA=-150,angleB=90,linecolor=darkred](1.5,1.366)(.5,0)
\pccurve[angleA=90,angleB=-150,linecolor=darkred](1.5,0)(2,.5)
\psline[linecolor=emgreen,linewidth=2pt](.25,3.531)(.5,3.098)
\psline[linewidth=2pt](.5,3.098)(1,2.232)
\psline[linecolor=emgreen,linewidth=2pt](1,2.232)(1.5,1.366)
\psline[linewidth=2pt](1.5,1.366)(2,.5)
\psline[linecolor=emgreen,linewidth=2pt](2,.5)(2.25,.067)
\psline[linecolor=emgreen,linewidth=2pt](-.25,3.531)(-.5,3.098)
\psline[linewidth=2pt](-.5,3.098)(-1,2.232)
\psline[linecolor=emgreen,linewidth=2pt](-1,2.232)(-1.5,1.366)
\psline[linewidth=2pt](-1.5,1.366)(-2,.5)
\psline[linecolor=emgreen,linewidth=2pt](-2,.5)(-2.25,.067)
\psline[linecolor=emgreen,linewidth=2pt](-1.75,0)(-1.5,0)
\psline[linewidth=2pt](-1.5,0)(-.5,0)
\psline[linecolor=emgreen,linewidth=2pt](-.5,0)(.5,0)
\psline[linewidth=2pt](.5,0)(1.5,0)
\psline[linecolor=emgreen,linewidth=2pt](1.5,0)(1.75,0)
\pscircle[linecolor=emgreen,linewidth=2.5pt](.25,3.531){.1}
\pscircle[linecolor=emgreen,linewidth=2.5pt](-.25,3.531){.1}
\pscircle[linecolor=emgreen,linewidth=2.5pt](2.25,.067){.1}
\pscircle[linecolor=emgreen,linewidth=2.5pt](-2.25,.067){.1}
\pscircle[linecolor=emgreen,linewidth=2.5pt](1.75,0){.1}
\pscircle[linecolor=emgreen,linewidth=2.5pt](-1.75,0){.1}
\endpspicture , \pspicture[.45](-2.5,-2)(2.5,1.5)
\pccurve[angleA=150,angleB=30,linecolor=darkred](1.366,1.366)(-1.366,1.366)
\pccurve[angleA=0,angleB=180,linecolor=darkred](-.5,.866)(.5,.866)
\pccurve[angleA=-120,angleB=150,linecolor=darkred](-1.866,.5)(-.5,-1.866)
\pccurve[angleA=120,angleB=-60,linecolor=darkred](-.5,-.866)(-1,0)
\pccurve[angleA=-120,angleB=60,linecolor=darkred](1,0)(.5,-.866)
\pccurve[angleA=30,angleB=-60,linecolor=darkred](.5,-1.866)(1.866,.5)
\psline[linewidth=1.5pt](.5,.866)(1.366,1.366)
\psline[linewidth=1.5pt](1.866,.5)(1,0)
\psline[linewidth=1.5pt](-.5,.866)(-1.366,1.366)
\psline[linewidth=1.5pt](-1.866,.5)(-1,0)
\psline[linewidth=1.5pt](-.5,-.866)(-.5,-1.866)
\psline[linewidth=1.5pt](.5,-1.866)(.5,-.866)
\psline[linecolor=emgreen,linewidth=1.5pt](-.5,.866)(-1,0)
\psline[linecolor=emgreen,linewidth=1.5pt](-1.866,.5)(-1.366,1.366)
\psline[linecolor=emgreen,linewidth=1.5pt](.5,.866)(1,0)
\psline[linecolor=emgreen,linewidth=1.5pt](1.866,.5)(1.366,1.366)
\psline[linecolor=emgreen,linewidth=1.5pt](-.5,-.866)(.5,-.866)
\psline[linecolor=emgreen,linewidth=1.5pt](.5,-1.866)(-.5,-1.866)
\pscircle[linecolor=emgreen,linewidth=2.5pt](1.616,.933){.1}
\pscircle[linecolor=emgreen,linewidth=2.5pt](-1.616,.933){.1}
\pscircle[linecolor=emgreen,linewidth=2.5pt](0,-1.866){.1}
\endpspicture
$$
\caption{A normal chord diagram and its corresponding polygonal
decomposition in Figure~\ref{decompo}.} \label{normal}
\end{figure}

All circular prime webs of a fixed $N$ tuple of even integers $(a_1,
a_2, \cdots, a_N)$, where $a_i\ge 4$, can be found the following
way. Because of $3$-connectivity of prime webs, for each polygon in
a circular polygonal decomposition, there exists a unique edge
between the polygon and the exterior face. Then we cut open each
polygon along this unique edge. Then we have a polygon of size $N$
with $a_i$ points marking on sides in cyclic order, let us call it a
\emph{plate} and $N$ is the \emph{size} of the plate. The edges of
the color which has not been used in polygonal decomposition give us
a chord diagram of the plate. Moreover, this chord diagram is called
\emph{normal} if there is no chord connecting points in the same
edge of the plate. Normal chord diagrams are well understood by the
representation theory of $\mathfrak{sl}(2)$. We will explain normal
chord diagrams by the size of the plate. If $N=2$, the only possible
normal chord diagram of the plate of size $2$ exists if $a_1$ is
equal to $a_2$. Thus if we fixed the number of the vertices of the
web, the prime web which has a polygonal decomposition of two
polygons is unique. For $N=3$, we say $(a_1,a_2,a_3)$ is an
\emph{admissible triple} if $|a_1-a_2|\le a_3\le a_1+a_2$. Let us
remind that all $a_i$ are even integers. For an admissible triple,
there exists a unique normal chord diagram. In the langauge of
representation theory \cite{kauffman, Turaevquantum}, let $V_{a_1}$
be an irreducible representation of $\mathfrak{sl}(2)$ of highest
weight $a_1$, then

$$\dim(\Inv( V_{a_1}\otimes V_{a_2}\otimes V_{a_3})) = \left\{
\begin{array}{cl} 1 & ~~\mathrm{if}~~ (a_1,a_2,a_3) ~~ \mathrm{is}
~\mathrm{an}~\mathrm{admissible}~\mathrm{triple}, \\
0 & ~~\mathrm{otherwise}. \end{array}\right.
$$

\begin{figure}
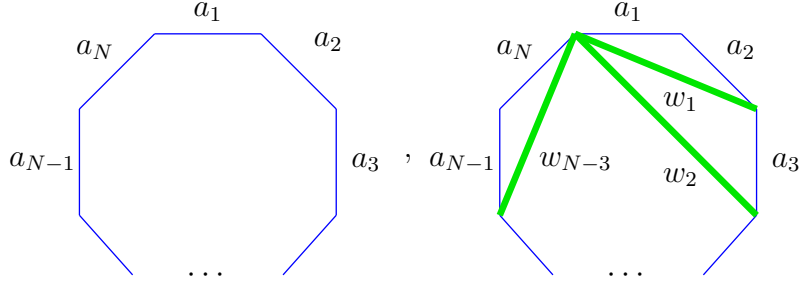

$$
\pspicture[.45](-2.4,-1.5)(2.6,2) \rput(0,2){\rnode{d1}{$a_1$}}
\rput(-1.5,1.5){\rnode{d2}{$a_N$}}
\rput(-2.2,0){\rnode{d3}{$a_{N-1}$}}
\rput(0,-1.5){\rnode{d4}{$\cdots$}} \rput(2.1,0){\rnode{d5}{$a_3$}}
\rput(1.6,1.6){\rnode{d6}{$a_2$}}
\qline(.707,1.707)(-.707,1.707)\qline(-1.707,.707)(-.707,1.707)
\qline(-1.707,.707)(-1.707,-.707)\qline(-1,-1.5)(-1.707,-.707)
\qline(1,-1.5)(1.707,-.707)\qline(1.707,-.707)(1.707,.707)
\qline(.707,1.707)(1.707,.707)
\endpspicture , \pspicture[.45](-2.8,-1.5)(2.4,2)
\rput(0,2){\rnode{d1}{$a_1$}} \rput(-1.5,1.5){\rnode{d2}{$a_N$}}
\rput(-2.2,0){\rnode{d3}{$a_{N-1}$}}
\rput(0,-1.5){\rnode{d4}{$\cdots$}} \rput(2.1,0){\rnode{d5}{$a_3$}}
\rput(1.5,1.5){\rnode{d6}{$a_2$}} \rput(.7,.8){\rnode{d5}{$w_1$}}
\rput(.7,-.2){\rnode{d5}{$w_2$}} \rput(-.7,0){\rnode{d5}{$w_{N-3}$}}
\qline(.707,1.707)(-.707,1.707)\qline(-1.707,.707)(-.707,1.707)
\qline(-1.707,.707)(-1.707,-.707)\qline(-1,-1.5)(-1.707,-.707)
\qline(1,-1.5)(1.707,-.707) \qline(1.707,-.707)(1.707,.707)
\qline(.707,1.707)(1.707,.707)
\psline[linecolor=emgreen,linewidth=2.5pt](-.707,1.707)(1.707,.707)
\psline[linecolor=emgreen,linewidth=2.5pt](-.707,1.707)(1.707,-.707)
\psline[linecolor=emgreen,linewidth=2.5pt](-.707,1.707)(-1.707,-.707)
\endpspicture
$$
\caption{A normal chord diagram interpolation of
equation~\ref{tensordecomp}.} \label{generaldecom}
\end{figure}

Let $\mathcal{C}$ be the set of all irreducible representations of
$\mathfrak{sl}(2)$. For $N\ge 4$, we use the following fact,

\begin{align}\nonumber
\dim&(\Inv(V_{a_1} \otimes V_{a_2} \otimes \ldots \otimes V_{a_N}))
= \sum_{V_{w_1}\in \mathcal{C}} \ldots \sum_{V_{w_{N-3}}\in
\mathcal{C}} \dim(\Inv(V_{a_1} \otimes V_{a_2} \otimes V_{w_1}))\\
&\dim(\Inv(V_{w_1} \otimes V_{a_3} \otimes V_{w_2})) \ldots
\dim(\Inv(V_{w_{N-3}} \otimes V_{a_{N-1}} \otimes
V_{a_N})).\label{tensordecomp}
\end{align}

Geometrically, this can be interpolated that we can count all normal
chord diagrams of fixed boundary as a sum of products of the normal
chord diagrams of triangles where the sums run on $\mathcal{C}$ by
the number of chords passing the line which represents an
irreducible representation $V_{w_i}$ as in Figure
\ref{generaldecom}. Even though $\mathcal{C}$ has infinitely many
elements, the actual sums run only on finitely many terms because
there are only finitely many admissible triples once we fix two
entries of the triples by the Clebsch-Gordan theorem
\cite{FultonHarrisgtm}.

\begin{figure}
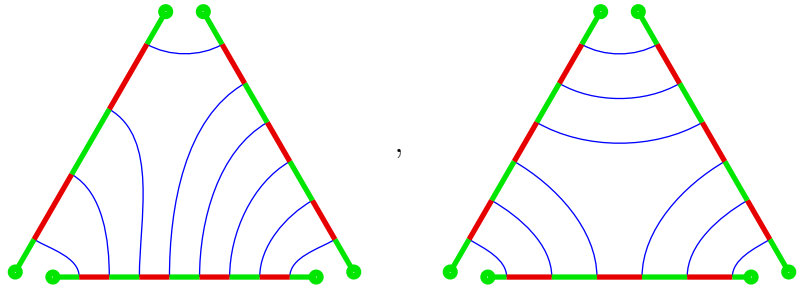

$$
\pspicture[.45](-2.8,0)(2.8,3.7)
\pccurve[angleA=-150,angleB=-30](.5,3.095)(-.5,3.098)
\pccurve[angleA=90,angleB=-30](-1.4,0)(-2,.5)
\pccurve[angleA=90,angleB=-30](-1,0)(-1.5,1.366)
\pccurve[angleA=90,angleB=-30](-.6,0)(-1,2.232)
\pccurve[angleA=90,angleB=-150](-.2,0)(.8,2.576)
\pccurve[angleA=90,angleB=-150](.2,0)(1.1,2.057)
\pccurve[angleA=90,angleB=-150](.6,0)(1.4,1.538)
\pccurve[angleA=90,angleB=-150](1,0)(1.7,1.019)
\pccurve[angleA=90,angleB=-150](1.4,0)(2,.5)
\psline[linecolor=emgreen,linewidth=2pt](.25,3.531)(.5,3.095)
\psline[linewidth=2pt,linecolor=darkred](2,.5)(1.7,1.019)
\psline[linecolor=emgreen,linewidth=2pt](1.7,1.019)(1.4,1.538)
\psline[linewidth=2pt,linecolor=darkred](1.4,1.538)(1.1,2.057)
\psline[linecolor=emgreen,linewidth=2pt](1.1,2.057)(.8,2.576)
\psline[linewidth=2pt,linecolor=darkred](.8,2.576)(.5,3.095)
\psline[linecolor=emgreen,linewidth=2pt](2,.5)(2.25,.067)
\psline[linecolor=emgreen,linewidth=2pt](-.25,3.531)(-.5,3.098)
\psline[linewidth=2pt,linecolor=darkred](-.5,3.098)(-1,2.232)
\psline[linecolor=emgreen,linewidth=2pt](-1,2.232)(-1.5,1.366)
\psline[linewidth=2pt,linecolor=darkred](-1.5,1.366)(-2,.5)
\psline[linecolor=emgreen,linewidth=2pt](-2,.5)(-2.25,.067)
\psline[linecolor=emgreen,linewidth=2pt](-1.75,0)(-1.4,0)
\psline[linewidth=2pt,linecolor=darkred](-1.4,0)(-1,0)
\psline[linecolor=emgreen,linewidth=2pt](-1,0)(-.6,0)
\psline[linewidth=2pt,linecolor=darkred](-.6,0)(-.2,0)
\psline[linecolor=emgreen,linewidth=2pt](-.2,0)(.2,0)
\psline[linewidth=2pt,linecolor=darkred](.2,0)(.6,0)
\psline[linecolor=emgreen,linewidth=2pt](.6,0)(1,0)
\psline[linewidth=2pt,linecolor=darkred](1,0)(1.4,0)
\psline[linecolor=emgreen,linewidth=2pt](1.4,0)(1.75,0)
\pscircle[linecolor=emgreen,linewidth=2.5pt](.25,3.531){.1}
\pscircle[linecolor=emgreen,linewidth=2.5pt](-.25,3.531){.1}
\pscircle[linecolor=emgreen,linewidth=2.5pt](2.25,.067){.1}
\pscircle[linecolor=emgreen,linewidth=2.5pt](-2.25,.067){.1}
\pscircle[linecolor=emgreen,linewidth=2.5pt](1.75,0){.1}
\pscircle[linecolor=emgreen,linewidth=2.5pt](-1.75,0){.1}
\endpspicture , \pspicture[.45](-2.8,0)(2.8,3.7)
\pccurve[angleA=-150,angleB=-30](.5,3.095)(-.5,3.095)
\pccurve[angleA=-150,angleB=-30](.8,2.576)(-.8,2.576)
\pccurve[angleA=-150,angleB=-30](1.1,2.057)(-1.1,2.057)
\pccurve[angleA=90,angleB=-30](-1.5,0)(-2,.5)
\pccurve[angleA=90,angleB=-30](-.9,0)(-1.7,1.019)
\pccurve[angleA=90,angleB=-30](-.3,0)(-1.4,1.538)
\pccurve[angleA=90,angleB=-150](.3,0)(1.4,1.538)
\pccurve[angleA=90,angleB=-150](.9,0)(1.7,1.019)
\pccurve[angleA=90,angleB=-150](1.5,0)(2,.5)
\psline[linecolor=emgreen,linewidth=2pt](.25,3.531)(.5,3.095)
\psline[linewidth=2pt,linecolor=darkred](2,.5)(1.7,1.019)
\psline[linecolor=emgreen,linewidth=2pt](1.7,1.019)(1.4,1.538)
\psline[linewidth=2pt,linecolor=darkred](1.4,1.538)(1.1,2.057)
\psline[linecolor=emgreen,linewidth=2pt](1.1,2.057)(.8,2.576)
\psline[linewidth=2pt,linecolor=darkred](.8,2.576)(.5,3.095)
\psline[linecolor=emgreen,linewidth=2pt](2,.5)(2.25,.067)
\psline[linecolor=emgreen,linewidth=2pt](-.25,3.531)(-.5,3.098)
\psline[linewidth=2pt,linecolor=darkred](-.5,3.095)(-.8,2.576)
\psline[linecolor=emgreen,linewidth=2pt](-.8,2.576)(-1.1,2.057)
\psline[linewidth=2pt,linecolor=darkred](-1.1,2.057)(-1.4,1.538)
\psline[linecolor=emgreen,linewidth=2pt](-1.4,1.538)(-1.7,1.019)
\psline[linewidth=2pt,linecolor=darkred](-1.7,1.019)(-2,.5)
\psline[linecolor=emgreen,linewidth=2pt](-2,.5)(-2.25,.067)
\psline[linecolor=emgreen,linewidth=2pt](-1.75,0)(-1.5,0)
\psline[linewidth=2pt,linecolor=darkred](-1.5,0)(-.9,0)
\psline[linecolor=emgreen,linewidth=2pt](-.9,0)(-.3,0)
\psline[linewidth=2pt,linecolor=darkred](-.3,0)(.3,0)
\psline[linecolor=emgreen,linewidth=2pt](.3,0)(.9,0)
\psline[linewidth=2pt,linecolor=darkred](.9,0)(1.5,0)
\psline[linecolor=emgreen,linewidth=2pt](1.5,0)(1.75,0)
\pscircle[linecolor=emgreen,linewidth=2.5pt](.25,3.531){.1}
\pscircle[linecolor=emgreen,linewidth=2.5pt](-.25,3.531){.1}
\pscircle[linecolor=emgreen,linewidth=2.5pt](2.25,.067){.1}
\pscircle[linecolor=emgreen,linewidth=2.5pt](-2.25,.067){.1}
\pscircle[linecolor=emgreen,linewidth=2.5pt](1.75,0){.1}
\pscircle[linecolor=emgreen,linewidth=2.5pt](-1.75,0){.1}
\endpspicture
$$
\caption{Normal chord diagrams of $(4,6,8)$ and $(6,6,6)$.}
\label{normalof9}
\end{figure}

\begin{figure}
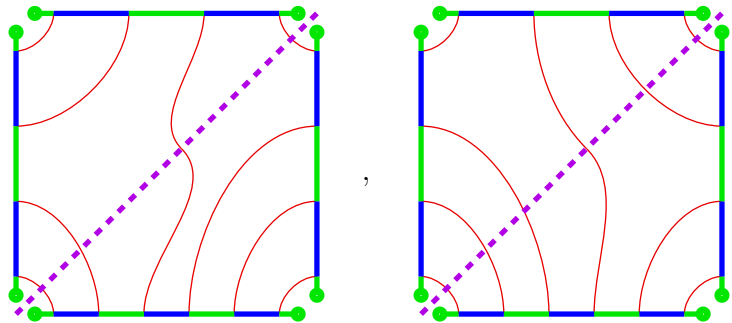

$$
\pspicture[.45](-2.2,-2.2)(2.6,2.2)\rput(1.5,-2){\rnode{a1}{$$}}
\pccurve[angleA=90,angleB=180,linecolor=darkred](1.5,-2)(2,-1.5)
\pccurve[angleA=90,angleB=180,linecolor=darkred](.9,-2)(2,-.5)
\pccurve[angleA=90,angleB=180,linecolor=darkred](.3,-2)(2,.5)
\pccurve[angleA=90,angleB=-45,linecolor=darkred](-.3,-2)(.2,.2)
\pccurve[angleA=135,angleB=-90,linecolor=darkred](.2,.2)(.5,2)
\pccurve[angleA=90,angleB=0,linecolor=darkred](-.9,-2)(-2,-.5)
\pccurve[angleA=90,angleB=0,linecolor=darkred](-1.5,-2)(-2,-1.5)
\pccurve[angleA=-90,angleB=0,linecolor=darkred](-.5,2)(-2,.5)
\pccurve[angleA=-90,angleB=0,linecolor=darkred](-1.5,2)(-2,1.5)
\pccurve[angleA=-90,angleB=180,linecolor=darkred](1.5,2)(2,1.5)
\psline[linecolor=emgreen,linewidth=2pt](-1.75,-2)(-1.5,-2)
\psline[linewidth=2pt](1.5,-2)(.9,-2)
\psline[linecolor=emgreen,linewidth=2pt](.9,-2)(.3,-2)
\psline[linewidth=2pt](.3,-2)(-.3,-2)
\psline[linecolor=emgreen,linewidth=2pt](-.3,-2)(-.9,-2)
\psline[linewidth=2pt](-.9,-2)(-1.5,-2)
\psline[linecolor=emgreen,linewidth=2pt](1.5,-2)(1.75,-2)
\psline[linecolor=emgreen,linewidth=2pt](2,-1.75)(2,-1.5)
\psline[linewidth=2pt](2,-1.5)(2,-.5)
\psline[linecolor=emgreen,linewidth=2pt](2,-.5)(2,.5)
\psline[linewidth=2pt](2,.5)(2,1.5)
\psline[linecolor=emgreen,linewidth=2pt](2,1.5)(2,1.75)
\psline[linecolor=emgreen,linewidth=2pt](-1.75,2)(-1.5,2)
\psline[linewidth=2pt](-1.5,2)(-.5,2)
\psline[linecolor=emgreen,linewidth=2pt](-.5,2)(.5,2)
\psline[linewidth=2pt](.5,2)(1.5,2)
\psline[linecolor=emgreen,linewidth=2pt](1.5,2)(1.75,2)
\psline[linecolor=emgreen,linewidth=2pt](-2,-1.75)(-2,-1.5)
\psline[linewidth=2pt](-2,-1.5)(-2,-.5)
\psline[linecolor=emgreen,linewidth=2pt](-2,-.5)(-2,.5)
\psline[linewidth=2pt](-2,.5)(-2,1.5)
\psline[linecolor=pup,linewidth=2pt,linewidth=2pt,linestyle=dashed](-2,-2)(2,2)
\psline[linecolor=emgreen,linewidth=2pt](-2,1.5)(-2,1.75)
\pscircle[linecolor=emgreen,linewidth=2.5pt](-2,1.75){.1}
\pscircle[linecolor=emgreen,linewidth=2.5pt](-2,-1.75){.1}
\pscircle[linecolor=emgreen,linewidth=2.5pt](2,1.75){.1}
\pscircle[linecolor=emgreen,linewidth=2.5pt](2,-1.75){.1}
\pscircle[linecolor=emgreen,linewidth=2.5pt](1.75,2){.1}
\pscircle[linecolor=emgreen,linewidth=2.5pt](-1.75,2){.1}
\pscircle[linecolor=emgreen,linewidth=2.5pt](1.75,-2){.1}
\pscircle[linecolor=emgreen,linewidth=2.5pt](-1.75,-2){.1}
\endpspicture , \pspicture[.45](-2.6,-2.2)(2.2,2.2)
\pccurve[angleA=90,angleB=180,linecolor=darkred](1.5,-2)(2,-1.5)
\pccurve[angleA=90,angleB=180,linecolor=darkred](.9,-2)(2,-.5)
\pccurve[angleA=90,angleB=-45,linecolor=darkred](.3,-2)(.2,.2)
\pccurve[angleA=135,angleB=-90,linecolor=darkred](.2,.2)(-.5,2)
\pccurve[angleA=90,angleB=0,linecolor=darkred](-.3,-2)(-2,.5)
\pccurve[angleA=90,angleB=0,linecolor=darkred](-.9,-2)(-2,-.5)
\pccurve[angleA=90,angleB=0,linecolor=darkred](-1.5,-2)(-2,-1.5)
\pccurve[angleA=-90,angleB=0,linecolor=darkred](-1.5,2)(-2,1.5)
\pccurve[angleA=-90,angleB=180,linecolor=darkred](1.5,2)(2,1.5)
\pccurve[angleA=-90,angleB=180,linecolor=darkred](.5,2)(2,.5)
\psline[linecolor=emgreen,linewidth=2pt](-1.75,-2)(-1.5,-2)
\psline[linewidth=2pt](1.5,-2)(.9,-2)
\psline[linecolor=emgreen,linewidth=2pt](.9,-2)(.3,-2)
\psline[linewidth=2pt](.3,-2)(-.3,-2)
\psline[linecolor=emgreen,linewidth=2pt](-.3,-2)(-.9,-2)
\psline[linewidth=2pt](-.9,-2)(-1.5,-2)
\psline[linecolor=emgreen,linewidth=2pt](1.5,-2)(1.75,-2)
\psline[linecolor=emgreen,linewidth=2pt](2,-1.75)(2,-1.5)
\psline[linewidth=2pt](2,-1.5)(2,-.5)
\psline[linecolor=emgreen,linewidth=2pt](2,-.5)(2,.5)
\psline[linewidth=2pt](2,.5)(2,1.5)
\psline[linecolor=emgreen,linewidth=2pt](2,1.5)(2,1.75)
\psline[linecolor=emgreen,linewidth=2pt](-1.75,2)(-1.5,2)
\psline[linewidth=2pt](-1.5,2)(-.5,2)
\psline[linecolor=emgreen,linewidth=2pt](-.5,2)(.5,2)
\psline[linewidth=2pt](.5,2)(1.5,2)
\psline[linecolor=emgreen,linewidth=2pt](1.5,2)(1.75,2)
\psline[linecolor=emgreen,linewidth=2pt](-2,-1.75)(-2,-1.5)
\psline[linewidth=2pt](-2,-1.5)(-2,-.5)
\psline[linecolor=emgreen,linewidth=2pt](-2,-.5)(-2,.5)
\psline[linewidth=2pt](-2,.5)(-2,1.5)
\psline[linecolor=emgreen,linewidth=2pt](-2,1.5)(-2,1.75)
\psline[linecolor=pup,linewidth=2pt,linewidth=2pt,linestyle=dashed](-2,-2)(2,2)
\pscircle[linecolor=emgreen,linewidth=2.5pt](-2,1.75){.1}
\pscircle[linecolor=emgreen,linewidth=2.5pt](-2,-1.75){.1}
\pscircle[linecolor=emgreen,linewidth=2.5pt](2,1.75){.1}
\pscircle[linecolor=emgreen,linewidth=2.5pt](2,-1.75){.1}
\pscircle[linecolor=emgreen,linewidth=2.5pt](1.75,2){.1}
\pscircle[linecolor=emgreen,linewidth=2.5pt](-1.75,2){.1}
\pscircle[linecolor=emgreen,linewidth=2.5pt](1.75,-2){.1}
\pscircle[linecolor=emgreen,linewidth=2.5pt](-1.75,-2){.1}
\endpspicture
$$
\caption{Normal chord diagrams of $(4,4,4)\cdot(4,4,6)$ and
$(4,4,6)\cdot(6,4,6)$.} \label{normalof9-1}
\end{figure}

\begin{figure}
\centerline{ \psfig{figure=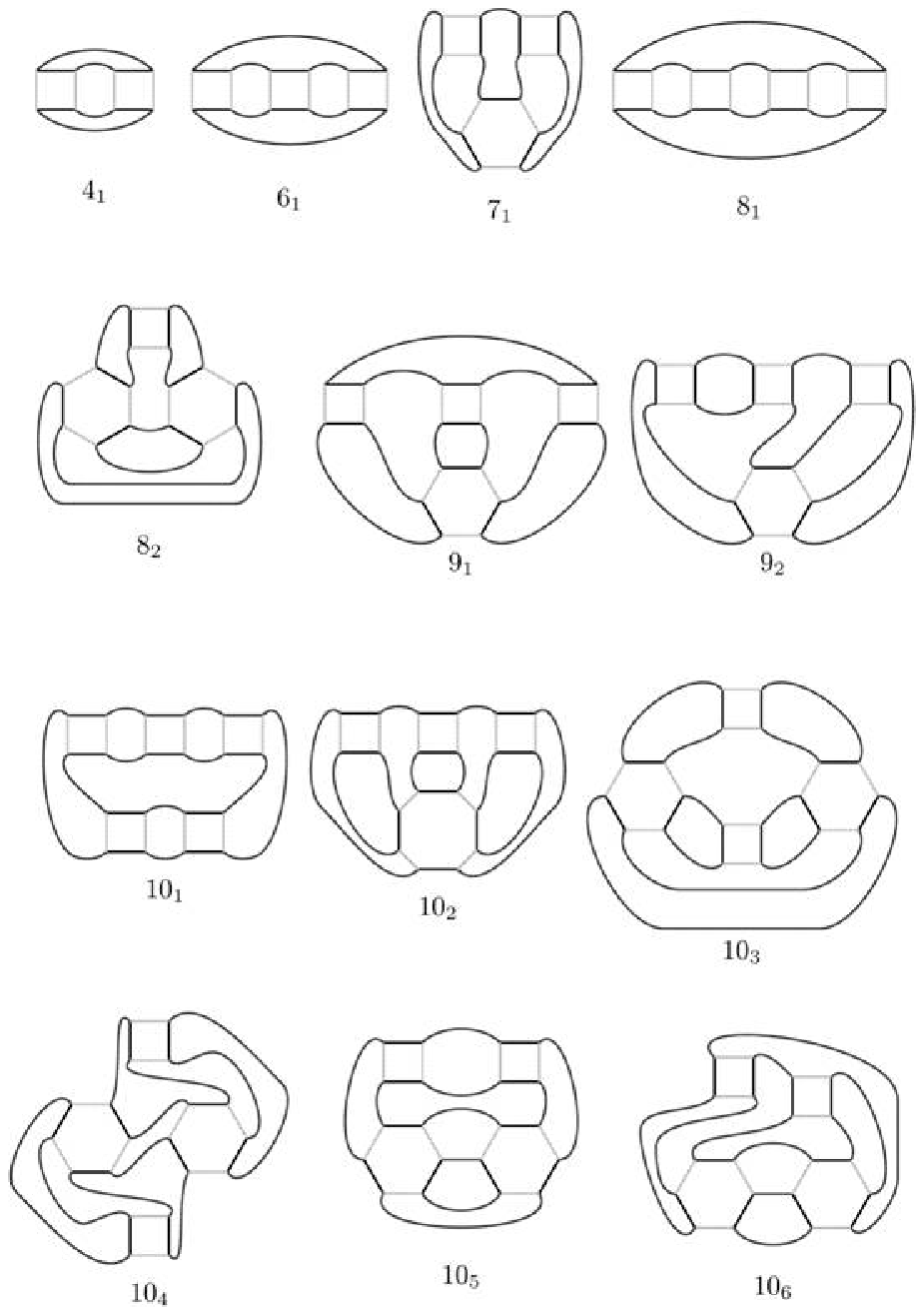,width=5in} } \vskip.1cm
\centerline{ \psfig{figure=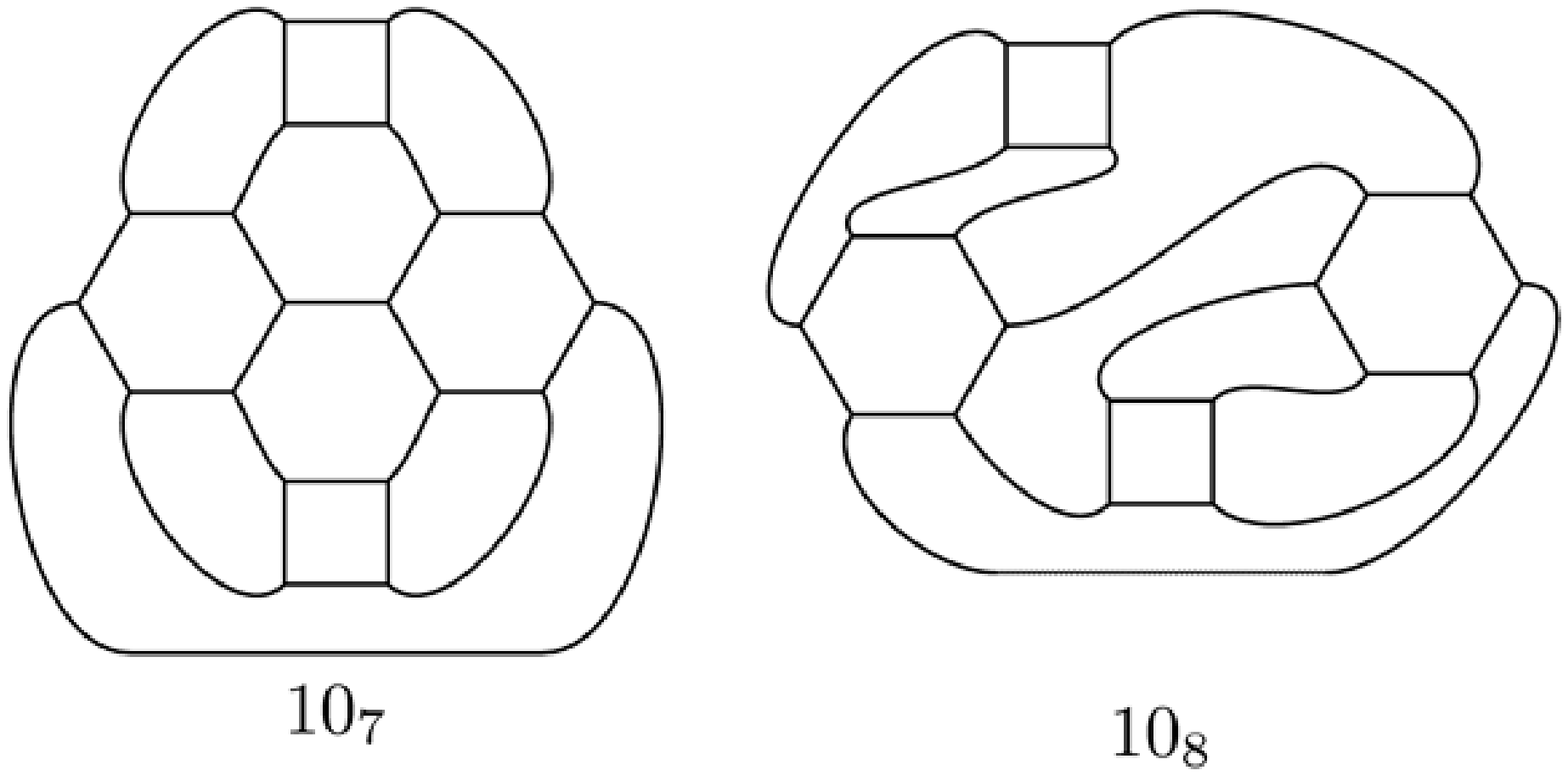,width=3in} } \caption{All prime
cubic bipartite planar graphs up to $20$ vertices.} \label{primes}
\end{figure}

For example, let us find all circular webs of $18$ vertices. Since
the sizes of every polygons are even and bigger than or equal to
$4$, we first find all even integral partitions of a given number of
vertices which satisfy these assumptions. For $18$, there are three
possible plate sizes, $2, 3$ and $4$. But $N=2$ can not be happened
because $\frac{18}{2}$ is not even. For $N=3$, we have three even
integral partitions of $18$ of length $3$, $(4,4,10), (4,6,8)$ and
$(6,6,6)$. But we can exclude $(4,4,10)$ because there does not
exist a normal chord diagrams of a triangle with $(4,4,10)$ markings
on the edges, or one can check $\dim (\Inv (V_4\otimes V_4\otimes
V_{10}))=0$. We find the chord diagrams of $(4,6,8)$ and $(6,6,6)$
in Figure \ref{normalof9}. The chord diagram of $(6,6,6)$ makes the
prime web $9_1$ and the chord diagram of $(4,6,8)$ makes the prime
web $9_2$ in Figure \ref{primes}. For $N=4$, the only even integral
partition of $18$ of length $4$ is $(4,4,4,6)$. From equation
\ref{tensordecomp}, we find that there are four possible nonzero
cases $(4,4,4)\cdot(2,4,6)$, $(4,4,4)\cdot(4,4,6)$,
$(4,4,6)\cdot(6,4,6)$ and $(4,4,8)\cdot(8,4,6)$. But one can easily
see that the first and the last are not prime. For the other we find
the chord diagrams of $(4,4,4)\cdot(4,4,6)$ and
$(4,4,6)\cdot(6,4,6)$ in Figure \ref{normalof9-1}. Both chord
diagrams make the prime web $9_2$ in Figure \ref{primes}.

\begin{figure}
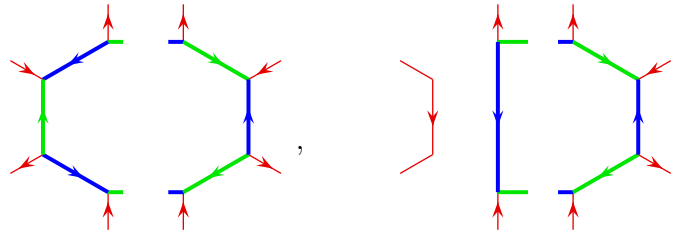

$$
\pspicture[.45](-2.5,-2)(2.5,1.5)
\pcline[linecolor=darkred](2.299,.75)(1.866,.5)\mda
\pcline[linecolor=darkred](1,1)(1,1.5)\mda
\pcline[linecolor=darkred](0,1)(0,1.5)\mda
\pcline[linecolor=darkred](-1.299,.75)(-.866,.5)\mda
\pcline[linecolor=darkred](-.866,-.5)(-1.299,-.75)\mda
\pcline[linecolor=darkred](0,-1.5)(0,-1)\mda
\pcline[linecolor=darkred](1,-1.5)(1,-1)\mda
\pcline[linecolor=darkred](1.866,-.5)(2.299,-.75)\mda
\psline[linecolor=emgreen,linewidth=1.5pt](0,1)(.2,1)
\psline[linecolor=emgreen,linewidth=1.5pt](0,-1)(.2,-1)
\psline[linewidth=1.5pt](1,1)(.8,1)
\psline[linewidth=1.5pt](1,-1)(.8,-1)
\pcline[linecolor=emgreen,linewidth=1.5pt](1,1)(1.866,.5)\mdb
\pcline[linewidth=1.5pt](0,1)(-.866,.5)\middlearrow
\pcline[linewidth=1.5pt](-.866,-.5)(0,-1)\middlearrow
\pcline[linecolor=emgreen,linewidth=1.5pt](1.866,-.5)(1,-1)\mdb
\pcline[linewidth=1.5pt](1.866,-.5)(1.866,.5)\middlearrow
\pcline[linecolor=emgreen,linewidth=1.5pt](-.866,-.5)(-.866,.5)\mdb
\endpspicture , \pspicture[.45](-2.5,-2)(2.5,1.5)
\pcline[linecolor=darkred](2.299,.75)(1.866,.5)\mda
\pcline[linecolor=darkred](1,1)(1,1.5)\mda
\pcline[linecolor=darkred](0,1)(0,1.5)\mda
\pcline[linecolor=darkred](-1.299,.75)(-.866,.5)
\pcline[linecolor=darkred](-.866,.5)(-.866,-.5)\mda
\psline[linecolor=darkred](-.866,-.5)(-1.299,-.75)
\pcline[linecolor=darkred](0,-1.5)(0,-1)\mda
\pcline[linecolor=darkred](1,-1.5)(1,-1)\mda
\pcline[linecolor=darkred](1.866,-.5)(2.299,-.75)\mda
\psline[linecolor=emgreen,linewidth=1.5pt](0,1)(.4,1)
\psline[linecolor=emgreen,linewidth=1.5pt](0,-1)(.4,-1)
\psline[linewidth=1.5pt](1,1)(.8,1)
\psline[linewidth=1.5pt](1,-1)(.8,-1)
\pcline[linecolor=emgreen,linewidth=1.5pt](1,1)(1.866,.5)\mdb
\pcline[linewidth=1.5pt](0,1)(0,-1)\middlearrow
\pcline[linecolor=emgreen,linewidth=1.5pt](1.866,-.5)(1,-1)\mdb
\pcline[linewidth=1.5pt](1.866,-.5)(1.866,.5)\middlearrow
\endpspicture
$$
\caption{A pushing move.} \label{pushing}
\end{figure}

\begin{figure}
$$
\pspicture[.45](-2.5,-3.2)(3,1.2)
\rput(-1.75,-1){\rnode{d1}{$\ldots$}}
\rput(0,-2.8){\rnode{d1}{$\mathrm{level}~2$}}
\rput(0,,8){\rnode{d1}{$\mathrm{exterior}$}}
\psline[linecolor=darkred](2,1)(1.5,.5)
\psline[linecolor=darkred](-1.5,.5)(-2,1)
\psline[linecolor=darkred](-2,0)(-2.5,0)
\psline[linecolor=darkred](2.5,0)(2,0)
\psline[linecolor=darkred](-1.5,-2.5)(-1,-2)
\psline[linecolor=darkred](1.5,-2.5)(1,-2)
\psline[linecolor=emgreen, linewidth=1.5pt](1.5,.5)(-1.5,.5)
\psline[linewidth=1.5pt](-1.5,.5)(-2,0)
\psline[linecolor=emgreen,linewidth=1.5pt](-2,0)(-2,-.5)
\psline[linecolor=emgreen,linewidth=1.5pt](-1.5,-1.5)(-1,-2)
\psline[linecolor=emgreen,linewidth=1.5pt](1,-2)(1.5,-1.5)
\psline[linecolor=emgreen,linewidth=1.5pt](2,-.5)(2,0)
\psline[linewidth=1.5pt](2,0)(1.5,.5)
\psline[linecolor=emgreen](-1.5,-2.5)(-2,-2.5)
\psline[linecolor=emgreen](1.5,-2.5)(2,-2.5)
\psline[linewidth=1.5pt](-1,-2)(1,-2) \psline(-1.5,-2.5)(-1.5,-3)
\psline(1.5,-2.5)(1.5,-3)
\endpspicture , \pspicture[.45](-3,-3.2)(2.5,1.2)
\rput(-1.75,-1){\rnode{d1}{$\ldots$}}
\rput(0,-2.8){\rnode{d2}{$\mathrm{level}~1$}}
\rput(0,,8){\rnode{d3}{$\mathrm{exterior}$}}
\rput(1.75,-1){\rnode{d4}{$\ldots$}}
\psline[linecolor=darkred](2,1)(1.5,.5)
\psline[linecolor=darkred](-1.5,.5)(-2,1)
\psline[linecolor=darkred](-2,0)(-2.5,0)
\psline[linecolor=darkred](2.5,0)(2,0)
\psline[linecolor=darkred](-1.5,-2.5)(-1,-2)
\psline[linecolor=darkred](1.5,-2.5)(1,-2)
\psline[linecolor=emgreen,linewidth=1.5pt](1.5,.5)(.5,.5)
\psline[linecolor=darkred](-.5,.5)(.5,.5)
\psline[linecolor=emgreen,linewidth=1.5pt](-1.5,.5)(-.5,.5)
\psline[linewidth=1.5pt](-1.5,.5)(-2,0)
\psline[linecolor=emgreen,linewidth=1.5pt](-2,0)(-2,-.5)
\psline[linecolor=emgreen,linewidth=1.5pt](-1.5,-1.5)(-1,-2)
\psline[linecolor=emgreen,linewidth=1.5pt](1,-2)(1.5,-1.5)
\psline[linecolor=emgreen,linewidth=1.5pt](2,-.5)(2,0)
\psline[linewidth=1.5pt](2,0)(1.5,.5)
\psline[linecolor=emgreen](-1.5,-2.5)(-2,-2.5)
\psline[linewidth=2.5pt](1,-2)(.5,-2)(.5,.5)
\psline[linewidth=2.5pt](-1,-2)(-.5,-2)(-.5,.5)
\psline(-1.5,-2.5)(-1.5,-3) \psline(1.5,-2.5)(1.5,-3)
\psline[linecolor=emgreen](1.5,-2.5)(2,-2.5)
\endpspicture
$$
\caption{A local shape around a polygon of level $2$ and a converse
of a pushing move.} \label{level2shape}
\end{figure}

To deal with non-circular webs, we define a \emph{pushing move} as
shown in Figure~\ref{pushing}. Then we find all non-circular webs
can be obtained by a finite sequence of pushing moves from
circular webs. First one can see any non-circular polygonal
decompositions of a web has a polygon of level $2$. The local
shape around a polygon of level $2$ is given in the left hand side
of Figure \ref{level2shape}. One can see that the there exists at
least one polygon of the same color type adjacent to the face
between the level $2$ polygon and the exterior polygon. Then if we
do a converse of a pushing move, then the sizes of new polygons,
$*$ in Figure \ref{level2shape}, are bigger than or equal to $4$.
Thus the resulting one is again a prime web. To finish the proof,
we induct on the number of polygons of level bigger than $1$ in a
non-circular polygonal decomposition of a non-circular prime web.
For a polygon of level $2$, we can see it can be obtained by a
pushing move from a polygonal decomposition with one less polygons
of level bigger than $1$. Inductively, it completes the proof.

For a fixed number $N$ of vertices, there is an upper bound $f(N)$
such that all prime webs of $N$ vertices can be found from
circular prime webs of $f(N)$ vertices. One of the easiest
estimation is $N+\lceil \frac{N}{8} \rceil$ because at each level
we must have at least two polygons. In fact, this can be improved
drastically because $f(20)$ is just $22$.

\subsection{Prime cubic bipartite
planar graphs up to $20$ vertices.}

Using the method described the above, we can find a list of all
prime webs up to $20$ vertices in Figure \ref{primes}. In Table
\ref{table1} and \ref{table2}, we list their quantum
$\mathfrak{sl}(3)$ invariants, three polygonal descriptions and
circularness of prime webs up to $20$ vertices. It is easy to find
all circular webs up to $22$ vertices, there are $8$ circular prime
webs of $22$ vertices, surprisingly these are all prime webs of $22$
vertices too. From these webs, we find all non-circular prime webs
up to $20$ vertices.

\section{Applications and Discussions}
\label{discussion}

\subsection{Symmetry of graphs and quantum $\mathfrak{sl}(3)$
invariants} \label{application}

A link $L$ in $S^3$ is $n$-\emph{periodic} if there exists a
periodic homeomorphism $h$ of order $n$ such that $fix(h)\cong S^1,
h(L)=L$ and $fix(h)\cap L=\emptyset$ where $fix(h)$ is the set of
fixed points of $h$. By a positive answer of the Smith conjecture,
if we consider $S^3$ as $\mathbb{R}^3\cup \{\infty\}$, we can assume
that $h$ is a rotation by $2\pi/n$ angle around the $z-$axis. If $L$
is a periodic link, we denote its quotient link by $L/h$. If a link
$L$ admits an orientation preserving action $\Gamma$ of order $n$,
then there are relations between the classical link polynomials such
as the Alexander polynomial and the Jones polynomial of link $L$ and
its quotient link $L/\Gamma$ \cite{mualexander, mujones}. For
quantum $\mathfrak{sl}(3)$ links invariants, it was shown
in~\cite{chbili}
$$P_L(q)\cong (P_{L/\Gamma}(q))^n \hskip ,7cm modulo \hskip .2cm
\mathcal{I}_n,$$ where $\mathcal{I}_n$ is the ideal of
$\mathbb{Z}[q^{\pm \frac 12}]$ generated by $n$ and $[3]^n-[3]$. It
has been generalized for the quantum $\mathfrak{sl}(n)$ link
invariants~\cite{JKsln}. In fact, the original invariants belongs to
$\mathbb{C}[q^{\pm \frac 12}]$ but later it was shown that it really
is a polynomial in
$\mathbb{Z}[q^{\frac12}+q^{-\frac12}]$~\cite{Leintegral}.

For planar cubic bipartite graphs, some of symmetries can be
orientation reversing. But the idea of the proof~\cite{chbili} still
works in general with one exception.  If the fundamental domain of
the action is a basis web with the given boundary, then there does
not exist any relation or expansion which can be applied repeatedly,
thus Theorem~\ref{symthm} might not work. In fact, we will
demonstrate this in Example \ref{ex1}.

 Let $G$ be a planar cubic bipartite
graph. Let $\Gamma$ be the group of orientation preserving
symmetries of $G$. Let $G/\Gamma$ be the quotient graph of $G$ by
$\Gamma$. A direct relation between $P_G(q)$ and $P_{G/\Gamma}(q)$
is very difficult to find in general. Thus we will look at this
relation modulo by an ideal in $\mathbb{Z}[q^{\pm \frac 12}]$.
Precisely we can state it as follows.

\begin{thm}
Let $G$ be a planar cubic bipartite graph with the group of
symmetries $\Gamma$ of order $n$. Let $\Gamma_d$ be a subgroup of
$\Gamma$ of order $d$ such that the fundamental domain of
$G/\Gamma_d$ is not a basis web with the given boundary. Then
$$P_{G}(q) \equiv (P_{G/\Gamma_d}(q))^d \hskip ,7cm modulo \hskip .2cm
\mathcal{I}_d,$$ where $\mathcal{I}_d$ is the ideal of
$\mathbb{Z}[q^{\pm \frac 12}]$ generated by $d$ and $[3]^d-[3]$.
\label{symthm}
\end{thm}

\begin{exmp}
We look at an example $6_1$. \label{ex1}
\end{exmp} There are twelve vertices but only six
of them have the same directions. Thus, its all possible symmetries
can have order either $6$, $3$ or $2$. For $n=3$, we find
\begin{align*}
[2]^4[3]+2[2]^2[3] &\equiv [3](([3]+[1])^2 + 2([3]+[1])) \equiv
[3]([3]^2 +
4[3] + 3[1]) \\
&\equiv [3]([3]^2+[3]) \equiv [3]^3+[3]^2 \equiv[3]+[3]^2 \\
&\equiv [3]([3]+[1]) \equiv [3]([3]^3+3[3]^2+3[3]+[1]) \\
&\equiv [3]([3]+[1])^3 \equiv [3]^3[2]^6 \equiv ([3][2]^2)^3
~~\mathrm{mod}~~\mathcal{I}_3
\end{align*}
In fact, it does have a symmetry of order three by a rotation along
a point in a hexagon. For $n=2$, we find $$[2]^4[3] +2[2]^2[3]
\equiv ([2]^2[3])^2 ~~\mathrm{mod}~~\mathcal{I}_2.$$ Also it does
have a symmetry of order two by a rotation along a point in a
rectangle. The quantum $\mathfrak{sl}(3)$ invariant of the quotient
of $6_1$ is $-[2][3]$ but it is very difficult to check whether its
6th power is congruent to $[2]^4[3] +2[2]^2[3]$ or not. By a help of
a machine~\cite{Baek}, we can see that there does not exist an
$\alpha \in \mathbb{Z}[q^{\pm \frac 12}]$ such that
$$(\alpha)^6 \equiv [2]^4[3]+2[2]^2[3]
~~\mathrm{mod}~~\mathcal{I}_6$$ even though there do exist a
symmetry of order $6$ for $6_1$.

\begin{table}
\begin{tabular}{|c|c|c|c|}\hline
$\begin{matrix} \mathrm{Prime} \\ \mathrm{Web} \end{matrix}$  &
Quantum ~~$\mathfrak{sl}(3)$ invariant & Polygonal descriptions & Circularness\\
\hline $4_1$ & $2[2]^2[3]$ &  $\begin{matrix} 4+4 \\ 4+4 \\
4+4 \end{matrix}$ & Yes \\  \hline
$6_1$ & $[2]^4[3]+2[2]^2[3]$ & $\begin{matrix} 4+4+4 \\ 4+4+4 \\
6+6 \end{matrix}$ & Yes  \\ \hline
$7_1$ & $-4[2]^3[3]$ &$\begin{matrix} 4+4+6 \\ 4+4+6 \\4+4+6 \end{matrix}$ & Yes  \\
\hline $8_1$ & $[2]^6[3] +[2]^4[3]+2[2]^2[3]$ & $\begin{matrix} 4+4+4+4\\
4+4+4+4 \\8+8
\end{matrix}$ & Yes \\  \hline
$8_2$ & $3[2]^4[3]+2[2]^2[3]$ & $\begin{matrix} 4+4+4+4\\ 4+6+6 \\
4+6+6 \end{matrix}$ & Yes \\ \hline
$9_1$ & $-[2]^5[3] -6[2]^3[3]$ & $\begin{matrix} 4+4+4+6\\ 4+4+4+6 \\
6+6+6 \end{matrix}$ & Yes \\ \hline
$9_2$ & $-2[2]^5[3] -4[2]^3[3]$ & $\begin{matrix} 4+4+4+6\\ 4+4+4+6 \\
4+6+8 \end{matrix}$ & Yes \\ \hline
$10_1$ & $[2]^8[3]+[2]^6+[2]^4[3]+2[2]^2[3]$ & $\begin{matrix} 4+4+4+4+4\\ 4+ 4+ 4 +4 +4\\
10+10 \end{matrix}$ & Yes \\ \hline
$10_2$ & $8[2]^4[3][3]$ & $\begin{matrix} 4+ 4+ 4 +8\\ 4+ 4+ 4 +8 \\
4+4+6+6 \end{matrix}$ & Yes \\ \hline
\end{tabular}
\vskip .5cm \caption{Quantum $\mathfrak{sl}(3)$ invariant, polygonal
descriptions and circularness of prime webs up to $20$ vertices.}
\label{table1}
\end{table}

\begin{table}
\begin{tabular}{|c|c|c|c|} \hline $\begin{matrix} \mathrm{Prime} \\ \mathrm{Web} \end{matrix}$  &
Quantum ~~$\mathfrak{sl}(3)$ invariant & Polygonal descriptions &
Circularness\\ \hline
$10_3$ & $[2]^6[3] +5[2]^4[3]+2[2]^2[3]$ & $\begin{matrix} 4+ 4+ 4 +4+4\\ 4+ 4+ 6 +6 \\
6+6+8 \end{matrix}$ & Yes \\
\hline $10_4$ & $8[2]^4[3]$ & $\begin{matrix} 4+ 4+ 4 +8\\ 4+ 4+ 6 +6 \\
4+4+6+6 \end{matrix}$ & Yes \\ \hline
$10_5$ & $6[2]^4[3][3] +3[2]^4[3]+2[2]^2[3]$ & $\begin{matrix} 4+ 4+ 4 +4+4\\ 4+ 4+ 6 +6 \\
4+8+8 \end{matrix}$ & Yes \\ \hline
$10_6$ & $7[2]^4[3] + 2[2]^2[3]$ & $\begin{matrix} 4+ 4+ 6 +6\\ 4+ 4+ 6 +6 \\
4+4+6+6 \end{matrix}$ & No \\ \hline
$10_7$ & $[2]^6[3]+5[2]^4[3]+2[2]^2[3]$ & $\begin{matrix} 4+ 4+ 6 +6\\ 4+ 4+ 6 +6 \\
4+4+6+6  \end{matrix}$ & No \\ \hline
$10_8$ & $8[2]^4[3]$ & $\begin{matrix} 4+ 4+ 6 +6\\ 4+ 4+ 6 +6 \\
4+4+6+6  \end{matrix}$ & No \\\hline
\end{tabular}
\vskip .5cm \caption{Table I continued.} \label{table2}
\end{table}

\subsection{Discussion}
\label{diss} The dual of a three-connected planar cubic graph is a
triangulation. In particular, the dual of prime webs are
three-connected planar Eulerian triangulations. The number of
these triangulations of a fixed vertices can be computed by a
C-program "plantri.c" \cite{plantri} and it has been computed up
to $68$ vertices \cite{Brinkmann}.

\begin{conj}
Let $k_n$ be the number of prime webs of $2n$ vertices. Then,
$$\lim_{n\rightarrow\infty} \frac{k_n}{k_{n-1}} = 3.829....$$
\end{conj}

One can ask what are the sufficient information to classify all
prime webs.

\begin{conj}
The quantum $\mathfrak{sl}(3)$ invariant, polygonal descriptions and
circularness completely determines all primes webs.
\end{conj}

For $n\ge 4$, there is the quantum $\mathfrak{sl}(n,\mathbb{C})$
invariants of webs, which are cubic weighted directed planar
graphs~\cite{JKsln}. It would be very interesting how one can
overcome the technical difficulty arose on the weights of edges.

\vskip .5cm \noindent{\bf Acknowledgements}

The authors would like to thank Greg Kuperberg for introducing the
subject and helpful discussion, Heungki Baek, Younghae Do, Myungsoo
Seo for their attention to this work. Also, the \TeX\, macro package
PSTricks~\cite{PSTricks} was essential for typesetting the equations
and figures.


\begin{thebibliography}{100}
\expandafter\ifx\csname
natexlab\endcsname\relax\def\natexlab#1{#1}\fi
\expandafter\ifx\csname url\endcsname\relax
  \def\url#1{\texttt{#1}}\fi
\expandafter\ifx\csname urlprefix\endcsname\relax\def\urlprefix{URL
}\fi

\bibitem{Baek} H. Baek, \textit{sl3-save.m}, Matlab-program.

\bibitem{barnette} D. Barnette, \textit{Conjecture 5, Recent progress in Combinatorics},
(Ed. W. T. Tutte), Academic Press. New York (1969), 343.

\bibitem{plantri} G. Brinkmann and B. McKay, \textit{Eplantri and fullgen}, C-program,
available at {\tt http://cs.anu.edu.au/\~{ }bdm/plantri/}.

\bibitem{Brinkmann} G. Brinkmann, private communication.

\bibitem{brown} W. Brown, \textit{Enumeration of triangulations of the disk},
Proc. London Math. Soc., 14 (1964), 746--768.

\bibitem{chbili} N. Chbili, \textit{The quantum SU(3) invariant of links and Murasugi's congruence},
Topology appl., 122 (2002), 479--485.


\bibitem{FultonHarrisgtm} W. Fulton and J. Harris, Representation theory,
Graduate Texts in Mathematics, 129, Springer-Verlag, New
York-Heidelberg-Berlin, 1991.

\bibitem{goodey} P. Goodey, \textit{Hamiltonian circuits in ploytopes with even sides},
Israel. J. Math. 22 (1975), 53--56.

\bibitem{GT1} J. Gross and T. Tucker, \textit{Topological graph theory},
Wiley-Interscience Series in discrete Mathematics and Optimization,
Wiley \& Sons, New York, 1987.

\bibitem{jaeger} F. Jaeger, \textit{A new invariant of plane bipartite cubic graphs},
Discrete math. 101 (1992), 149--164.

\bibitem{Jonespoly} V. Jones, \textit{Polynomial invariants of knots via von Neumann algebras},
Bull. Amer. Math. Soc., 12 (1985), 103--111.

\bibitem{JKsln} M. Jeong and D. Kim, \textit{Quantum $\mathfrak{sl}(n)$ link invariants}, preprint,
arXiv:math.GT/0506403.

\bibitem{kauffman} L. Kauffman and S. Lins, \textit{Temperley-Lieb
recoupling theory and invariants of 3-manifolds}, Annals of Math.
Studies, 134, Princeton University Press, Princeton, 1994.

\bibitem{Khovanovsl3} M. Khovanov, \textit{sl(3) link homology}, Algebr. Geom. Topol.,
4 (2004), 1045-1081.

\bibitem{Khovanovcolored} M. Khovanov, \textit{Categorifications of the colored Jones polynomial},
J. Knot Theory Ramifications, 14~(1) (2005) 111--130,
arXiv:math.QA/0302060.


\bibitem{Dongseokthesis} D. Kim, \textit{ Graphical Calculus on Representations of Quantum Lie
Algebras}, Thesis, UCDavis, 2003, arXiv:math.QA /0310143.

\bibitem{Kuperbergspiders} G. Kuperberg, \textit{Spiders for rank 2 {Lie} algebras}, Comm. Math. Phys.,
180(1) (1996), 109--151, arXiv:q-alg/9712003.

\bibitem{KKnotdual} M. Khovanov and G. Kuperberg, \textit{Web bases for $sl(3)$ are not dual canonical},
Pacific J. Math., 188(1) (1999), 129--153, arXiv:q-alg/9712046.

\bibitem{KL} J. Kim and J. Lee, \textit{Genus distributions for
bouquets of dipoles}, J. Korean Math. Soc., 35 (1998), 225--235.

\bibitem{Leintegral} T. Le, \textit{Integrality and symmetry of quantum Link invariants},
Duke Math. J., 102 (2000), 273--306.

\bibitem{Lickorishsu} W. Lickorish, \textit{Distinct 3-manifolds with all $SU(2)_q$ invariants the same},
Proc. Amer. Math. Soc., 117 (1993), 285--292.

\bibitem{mualexander} K. Murasugi, \textit{On periodic knots}, Comment. Math. Helv.,
46 (1971), 162--174.

\bibitem{mujones} K. Murasugi, \textit{The Jones polynomials of periodic links},
Pacific J. Math., 131 (1988), 319--329.

\bibitem{MOYHomfly} H. Murakami and T. Ohtsuki and S. Yamada, \textit{HOMFLY
polynomial via an invariant of colored plane graphs}, L'Enseignement
Mathematique, t., 44 (1998), 325--360.

\bibitem{OYquantum} T. Ohtsuki and S. Yamada, \textit{Quantum $su(3)$ invariants via linear skein theory},
J. Knot Theory Ramifications, 6(3) (1997), 373-404.

\bibitem{PSTricks}
T.~{Van Zandt}. PSTricks: {PostScript} macros for generic {\TeX},
available at {\tt ftp://ftp.princeton.edu/pub/tvz/}.

\bibitem{RT1} N. Reshetikhin and V. Turaev,
\textit{Ribbob graphs and their invariants derived from quantum
groups}, Comm. Math. Phys., 127 (1990), 1--26.

\bibitem{Sokolov} M. Sokolov,
\textit{On the absolute value of the SO(3)-invariant and other
summands of the Turaev-Viro invariant}, Banach Center Publ., 42,
Polish Acad. Sci. Knot theory (Warsaw, 1995), (1998), 395--408,
arXiv:q-alg/9601013.

\bibitem{Turaevquantum} V. Turaev, \textit{Quantum invariants of knots and 3-manifolds},
W. de Gruyter, Berlin, 1994.

\bibitem{Walshcounting} T. Walsh, \textit{Counting unlabelled three-connected and homeomorphically
irreducible two connected graphs}, J. Combin. Theory., 32 (1982),
12-32.

\bibitem{whitney} H. Whitney, \textit{A set of topological invariants for graphs},
Amer. J. Math., 55 (1933), 231-235.

\end{thebibliography}
\end{document}